\documentclass[final]{AeroBest_paper}

\usepackage{graphicx}
\usepackage{booktabs}
\usepackage{amssymb}

\usepackage{subcaption}

\usepackage{amsmath}
\usepackage{amsfonts}

\usepackage{hyperref} 
\hypersetup{pdfborder=0 0 0}

\usepackage{doi}

\usepackage[square,numbers,sort&compress]{natbib}
\title{THERMO-MECHANICAL LEVEL-SET TOPOLOGY OPTIMIZATION OF A LOAD CARRYING BATTERY PACK FOR ELECTRIC AIRCRAFT}

\author{Alexandre T.R. Guibert$^{\bf 1}$, Murtaza Bookwala$^{\bf 1}$, Ashley Cronk$^{\bf 2}$, Y. Shirley Meng$^{\bf 3}$, and H. Alicia Kim$^{\bf 1,2}$}

\heading{Alexandre T.R. Guibert, Murtaza Bookwala, Ashley Cronk, Y. Shirley Meng, and H. Alicia Kim}

\address{1: Structural Engineering Department\\
Jacobs School of Engineering\\
University of California San Diego\\
Matthews Lane, La Jolla, CA 92093, USA\\
\{aguibert,mbookwala,alicia\}@ucsd.edu,  http://m2do.ucsd.edu/\\
\
\\
2: Materials Science and Engineering Department \\
Jacobs School of Engineering\\
University of California San Diego\\
Matthews Lane, La Jolla, CA 92093, USA\\
acronk@eng.ucsd.edu, http://smeng.ucsd.edu/ \\
\
\\
3: Pritzker School of Molecular Engineering\\
The University of Chicago\\
5801 S Ellis Avenue, Chicago, 60637, Illinois, USA\\
shirleymeng@uchicago.edu, https://voices.uchicago.edu/smeng/}

\abstract{
A persistent challenge with the development of electric vertical take-off and landing vehicles (eVTOL) to meet flight power and energy demands is the mass of the load and thermal management systems for batteries. One possible strategy to overcome this problem is to employ optimization techniques to obtain a lightweight battery pack while satisfying structural and thermal requirements. In this work, a structural battery pack with high-energy-density cylindrical cells is optimized using the level-set topology optimization method. The heat generated by the batteries is predicted using a high-fidelity electrochemical model for a given eVTOL flight profile. The worst-case scenario for the battery's heat generation is then considered as a source term in the weakly coupled steady-state thermo-mechanical finite element model used for optimization. The objective of the optimization problem is to minimize the weighted sum of thermal compliance and structural compliance subjected to a volume constraint. The methodology is demonstrated with numerical examples for different sets of weights. The optimized results due to different weights are compared, discussed, and evaluated with thermal and structural performance indicators. The optimized pack topologies are subjected to a transient thermal finite element analysis to assess the battery pack's thermal response. 
}


\begin{document}

%

\section{Introduction}
\label{sec:intro}

The aerospace sector is undergoing a significant transformation due to the emergence of electric aircraft and the advancement of urban air mobility (UAM) concepts. One of the solutions to UAM is the electric Vertical Take-Off and Landing (eVTOL) vehicles, which have the potential to alleviate ground traffic congestion and provide a secure and environmentally friendly alternative to conventional individual transportation.
However, eVTOL vehicles pose fundamental engineering challenges such as noise control and power requirements. To make eVTOL vehicles a viable transportation option, there is a crucial need to enhance their power density while ensuring safety by carefully designing the vehicles' battery packs. By utilizing topology optimization, it is possible to optimize the distribution of materials in a battery pack with the objective of improving the structural and thermal performance while minimizing its weight which ultimately contributes to the overall efficiency of the aircraft. This can eliminate a mechanical and thermal load management system which is otherwise needed, thus reducing the mass and volume budget. Therefore, the integration of topology optimization in the development of eVTOL vehicles presents a promising opportunity to move closer to a sustainable and efficient form of UAM.

In this work, we present a methodology for the multiphysics optimization of a battery pack for electric aircraft using the level set method. The volumetric heat generation rate from the batteries is estimated using an electrochemical model and the worst-case scenario is considered for design optimization. The electrochemical model, heat transfer model, and elasticity model are introduced in Section \ref{electro}, \ref{heat}, and \ref{elasticity} respectively. The level-set topology optimization method and the optimization problem and workflow are discussed in Section \ref{TOSection}. The methodology is demonstrated with numerical examples in Section \ref{examples}.

\section{Multiphysics model}
\subsection{High fidelity electrochemical model} \label{electro}

The cell behavior is modeled using the Doyle Fuller Newman (DFN) model \cite{Doyle1993, Ai2019, Timms2021, Deshpande2012} which describes the electrochemical response of lithium-ion cells. It is solved for the electric potentials and lithium-ion concentration distribution in a cell given some electric loads such as an applied current, voltage, or power. The DFN model in its one-dimensional version assumes that the charge/discharge process is mainly unidirectional, i.e, only the effects from one current collector to another are considered. It also assumes that the solid particles are spherical. The lithium-ion diffusion in these particles is considered by having a pseudo second dimension associated with the particles' radius. The governing equations for the DFN model are the following: 
\begin{enumerate}
    \item Molar conservation in the electrolyte,
    \item Solid state Fickian diffusion in the the solid particles, 
    \item Charge conservation in the electrolyte and in the particles based on Ohm's law,
    \item Intercalation kinetics based on the Butler-Volmer equation which effectively couples the macro dimension and the pseudo dimension.
\end{enumerate}
The DFN model is implemented using the open-source library PyBaMM \cite{Sulzer2021, Andersson2019, Harris2020}.

\subsection{Heat transfer model} \label{heat}

The governing equation for the steady-state thermal model considering thermal conduction only is, 
\begin{equation}
\nabla \cdot (\kappa \nabla T) + Q = 0 \,,
\end{equation}
where $\kappa$ is the thermal conductivity, $T$ is the temperature, and $Q$ is the volumetric heat generation rate. The latter is determined by the electrochemical model by adding the volumetric heat generation due to Ohmic heating $Q_{Ohm}$ from the electrical resistance of the cell, irreversible heating from the electrochemical reactions $Q_{rxn}$, and reversible heating due to entropic changes $Q_{rev}$. The heat transfer model is implemented using the finite element method and automatic differentiation software, FEniCS and DOLFIN \cite{alnaes2015fenics, logg2010dolfin, logg2012automated}.

\subsection{Elasticity model} \label{elasticity}

The governing equation of the linear elasticity problem is,
\begin{equation}
 \nabla \cdot \boldsymbol \sigma(\boldsymbol u) + \boldsymbol b = 0 \,,
\end{equation}
where $\boldsymbol b$ is the body force, $\boldsymbol u$ is the displacement, and $\boldsymbol \sigma$ is the Cauchy stress tensor. The temperature distribution solved by the heat transfer model resulting in thermal strain $\varepsilon_T$,
\begin{equation}
    \varepsilon_T = \alpha \Delta T \,,
\end{equation}
where $\alpha$ is the coefficient of thermal expansion and $\Delta T$ is the temperature increase, that is, $\Delta T = T - T_{ref}$ where $T_{ref}$ is the reference temperature. The constitutive law then becomes, 
\begin{equation}
    \boldsymbol \sigma = \mathbb{C}:(\boldsymbol \varepsilon - \varepsilon_T\boldsymbol 1) \,,
\end{equation}
where $\mathbb{C}$ is the elastic constitutive tensor, $\boldsymbol 1$ is the unit tensor, and $\boldsymbol \varepsilon$ is the strain tensor defined as $\boldsymbol \varepsilon = \frac{1}{2}(\nabla\boldsymbol u   + \nabla\boldsymbol u^\top)$. The linear elasticity model is implemented in its weak form using FEniCS and DOLFIN \cite{alnaes2015fenics, logg2010dolfin, logg2012automated}.

\section{Topology optimization} \label{TOSection}

\subsection{Level-set topology optimization} \label{levelset}
The optimized design is described by an implicit function $\phi(\boldsymbol x)$ such that
\begin{equation}
    \begin{cases}
        \phi({\bf x}) \geq 0, {\bf x} \in \Omega \\
        \phi({\bf x}) = 0, {\bf x} \in \Gamma \\
        \phi({\bf x}) < 0, {\bf x} \not\in \Omega
    \end{cases}   
\end{equation}
where $\Omega$ is the solid domain such that $\Omega \subset \mathcal{D}$ where $\mathcal{D}$ is the design domain and $\mathcal{D}\setminus \Omega$ is the void domain. The boundary of the solid domain is $\Gamma$ described by the zero level-set.  The design is iteratively updated by the discretized Hamilton-Jacobi equation,
\begin{equation}
    \phi^{k+1}_i = \phi^k_i - \Delta t \vert \nabla\phi^k_i  \vert V_{n, i} \,,
\end{equation}
where $k$ is the iteration number, $i$ is a boundary point, $\Delta t$ is the pseudo time-step, and $V_n$ is the design velocity normal to $\Gamma$ \cite{dunning2015introducing, sivapuram2016simultaneous}. The finite element analysis is carried out on a fixed grid, i.e., a non-conforming mesh is used. Thus, as the design evolves, the zero level-set cuts finite elements such that for a given element some nodes of that element might be inside the optimized design and some might be outside. To address this issue, the following material interpolation scheme is used,
\begin{equation}
\begin{split}
    \kappa_e = \left(\gamma_{min}(1 - \gamma_e) + \gamma_{e}\right) \kappa \,, \\
    E_e = \left(\gamma_{min}(1 - \gamma_e) + \gamma_{e}\right) E \,,
\end{split}
\end{equation}
where $\gamma_e$ is the fraction of the element where $\phi({\bf x}) > 0$, and $\gamma_{min}$ is a small value, taken as $10^{-4}$ for this study, to ensure that the stiffness and conductivity matrices are non-singular. The optimization is implemented using ParaLeSTO \cite{kambampati_sandilya_2023_7613753} which is a modular open-source level-set topology optimization code  \cite{hyun2022development, guibert2023implementation}. 

\subsection{Optimization workflow and formulation} \label{workflow}

\begin{figure}[!htb]
      \centering
      \includegraphics[width=0.9\textwidth]{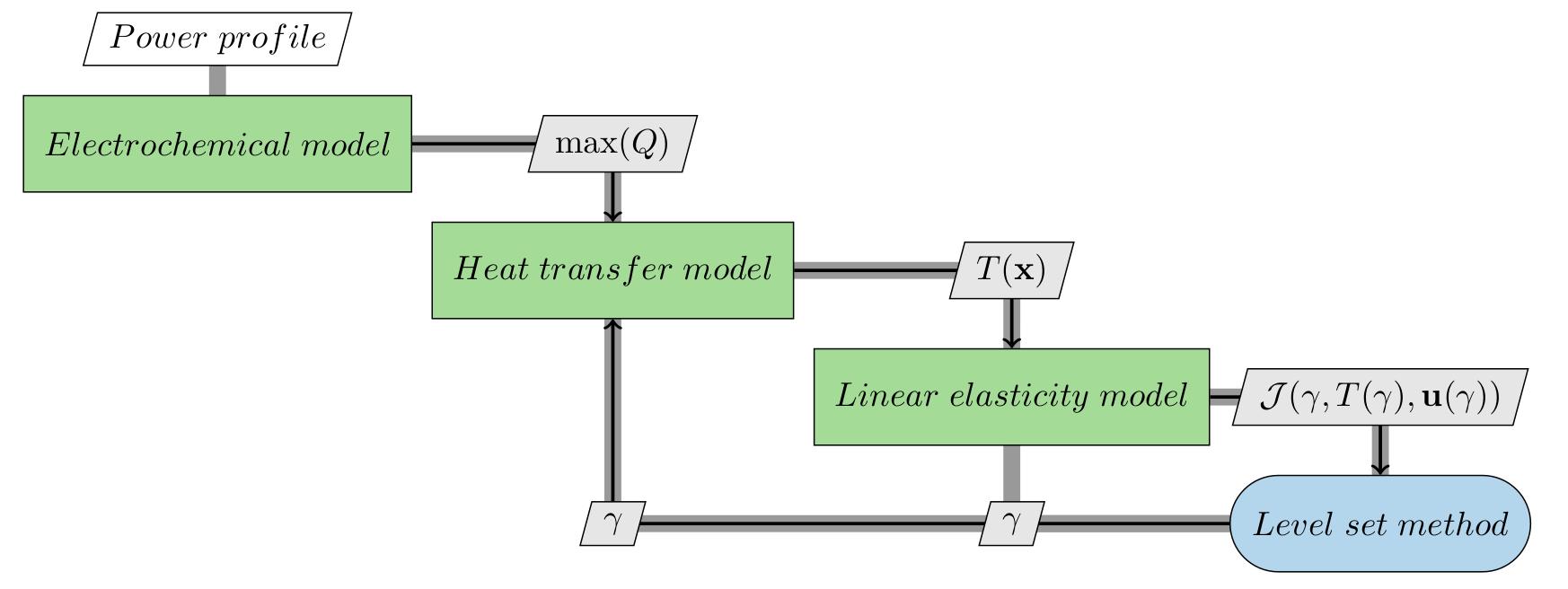}
   \caption{Optimization workflow}
   \label{fig:myXDSM}
\end{figure}

The objective function of the optimization problem $\mathcal{J}$ is the weighted sum of the normalized thermal and structural compliance, i.e., 
\begin{equation}
    \mathcal{J}({ \gamma},{T}({\gamma}), {\bf u}({\gamma}))  \triangleq k_1 \frac{C_S}{C_{S0}} + k_2 \frac{C_T}{C_{T0}} \,,
    \label{objective1}
\end{equation}
where $C_S$ and $C_T$ are the structural compliance and thermal compliance defined as
\begin{equation}
    \begin{split}
        C_S = \boldsymbol F_S^\top\boldsymbol U \,,\\
        C_T = \boldsymbol F_T^\top\boldsymbol T \,,\\
    \end{split}
\end{equation}
where $\boldsymbol F_T$ and $\boldsymbol F_S$ are the thermal force vector and structural force vector respectively. The normalization factors $C_{S0}$ and $C_{T0}$ are reinitialized every 5 iterations. Moreover, let $k \equiv k_1 = 1 - k_2$. Thus, Eq.(\ref{objective1}) becomes
\begin{equation}
    \mathcal{J}({ \gamma},{T}({\gamma}), {\bf u}({\gamma})) = k \frac{C_S}{C_{S0}} + (1-k) \frac{C_T}{C_{T0}} \,.
    \label{objective2}
\end{equation}
Additionally, the volume for the optimized part is constrained such that a fraction $\xi$ of the design domain is prescribed. Finally, the optimization problem to be solved can be expressed as
\begin{eqnarray}
  {\rm Minimize}   && \mathcal{J}({ \gamma},{T}({\gamma}), {\bf u}({\gamma}))            \nonumber           \\
  {\rm w.r.t.}     && {\gamma}                                     \label{eq:minimize} \\
  {\rm subject~to} &&    \operatorname{Vol}(\Omega) \leq \xi \operatorname{Vol}(\mathcal{D})  \nonumber    \\
                   &&    {\cal R}_1({\gamma},{T}({\gamma})) = 0   \nonumber \\
                   &&    {\cal R}_2({\gamma},{\bf u}({\gamma})) = 0   \nonumber
\end{eqnarray}
where ${\cal R}_1$ and ${\cal R}_2$ are the residuals of the thermal and elasticity models respectively. The optimization workflow is presented in Fig.\ref{fig:myXDSM}.

\section{Numerical results} \label{examples}

\subsection{Set up and boundary conditions}
\begin{figure}[!h]
\centering
\begin{subfigure}{0.49\textwidth}
    \includegraphics[width=\textwidth, trim={0 0 0 2.9cm}]{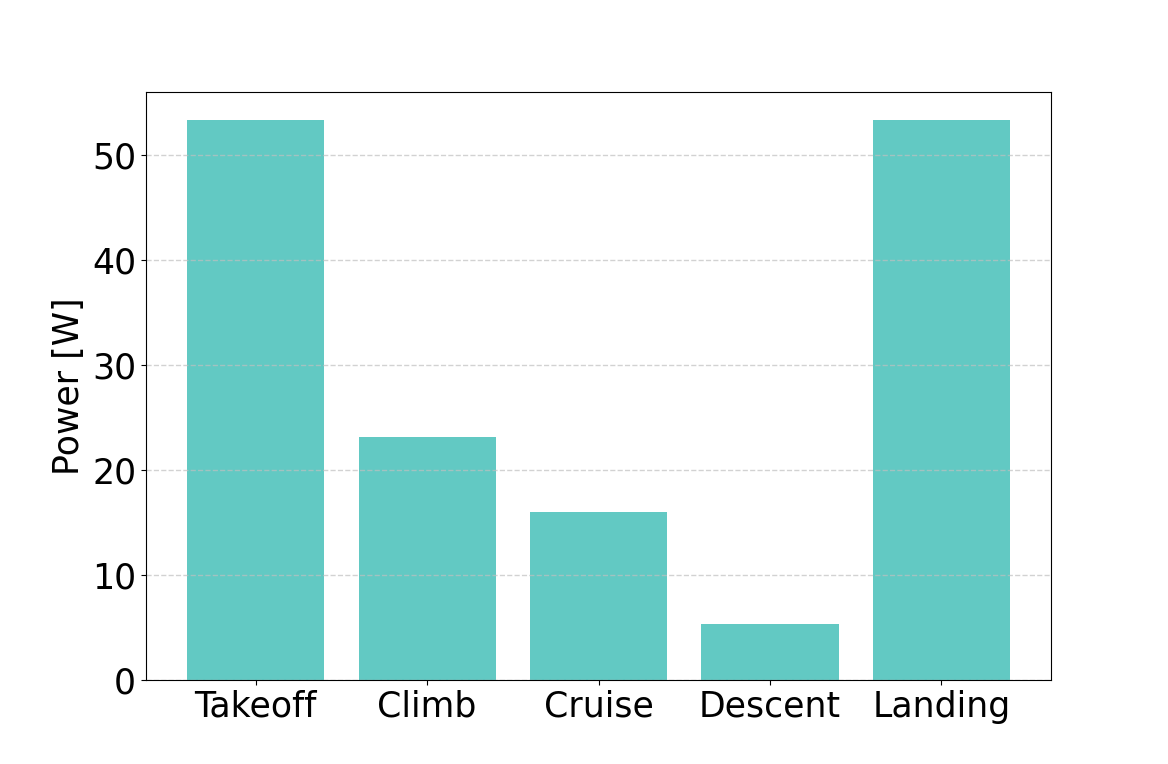}
    \caption{Power profile}
\end{subfigure}
\hfill
\begin{subfigure}{0.49\textwidth}
    \includegraphics[width=\textwidth, trim={0 0 0 2.9cm}]{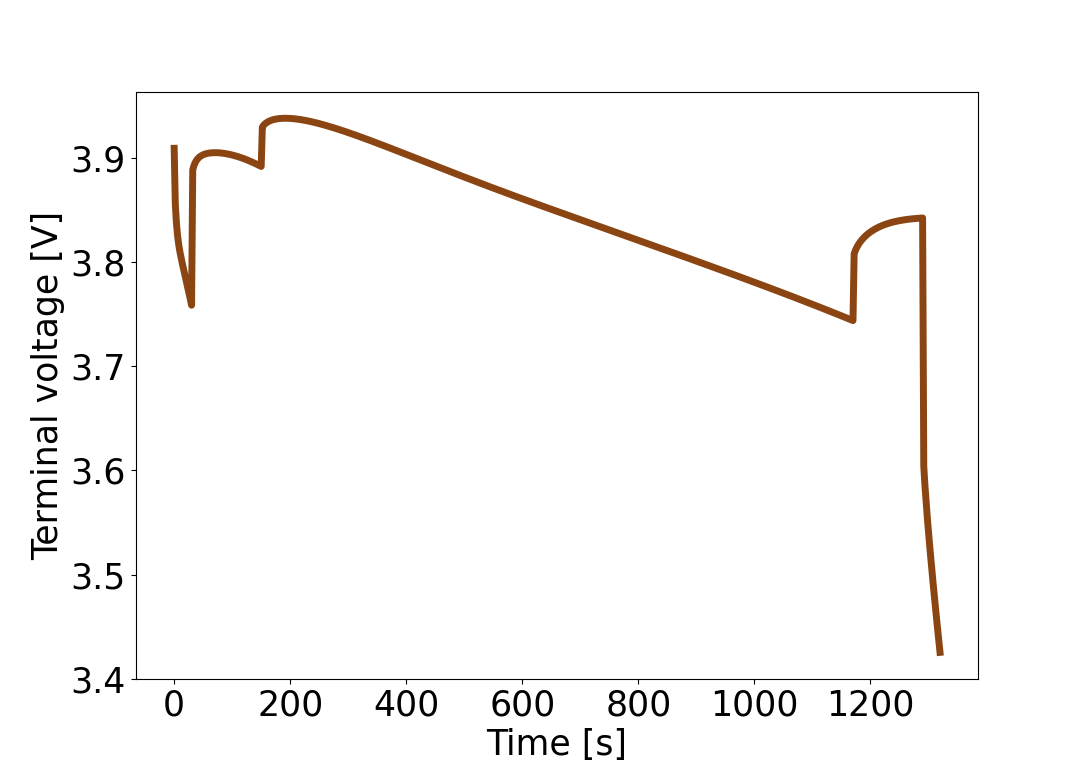}
    \caption{Voltage profile}
\end{subfigure}
\hfill
\begin{subfigure}{\textwidth}
    \centering
    \includegraphics[width=0.85\textwidth]{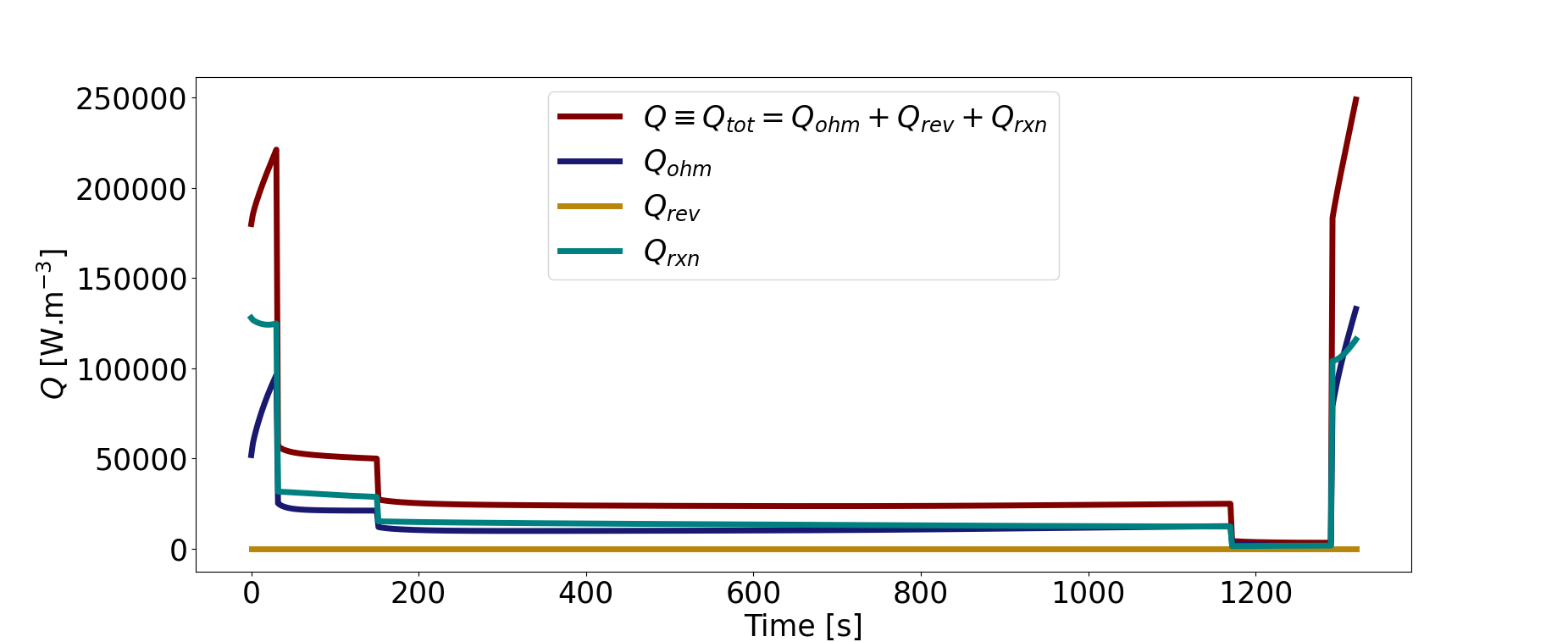}
    \caption{Heat generation from power profile}
\end{subfigure}
\caption{Electrochemical response of 21700 cells for a given power profile}
\label{fig:Electro}
\end{figure}

The LG CHEM INR21700-M50L cell is used for its high energy density and the properties are taken from \cite{Chen2020}.
The power profile is computed using \cite{yang2021challenges} where the energy is taken as 17.8 Wh \cite{batemo_2022}. The power profile, the voltage profile obtained with the electrochemical model, and the estimated volumetric heat generation rate are presented in Fig.\ref{fig:Electro}.

\begin{figure}[!h]
\centering
\begin{subfigure}{0.49\textwidth}
    \includegraphics[width=\textwidth, trim={6cm 0 6cm 0}]{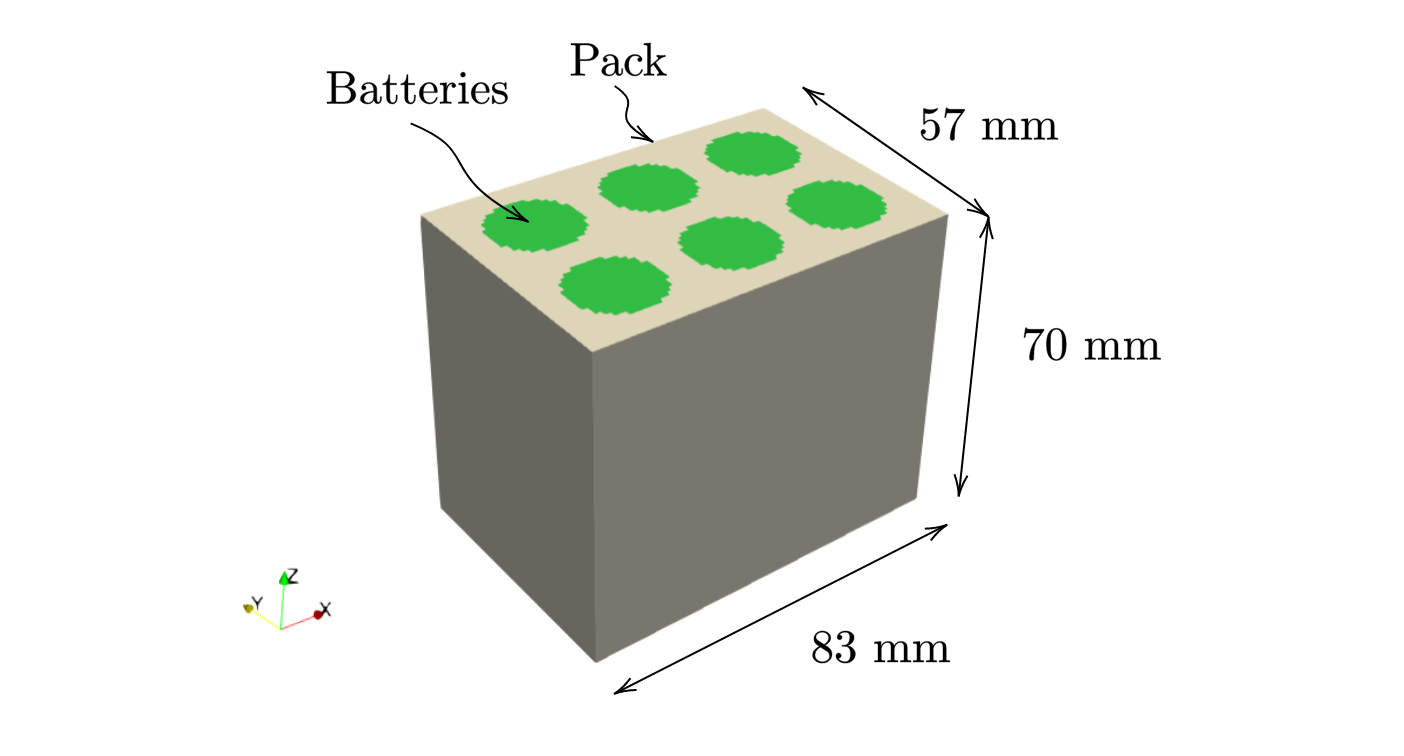}
    \caption{Geometry of the problem}
\end{subfigure}
\hfill
\begin{subfigure}{0.49\textwidth}
    \includegraphics[width=\textwidth, trim={4.5cm 0 4.5cm 0}]{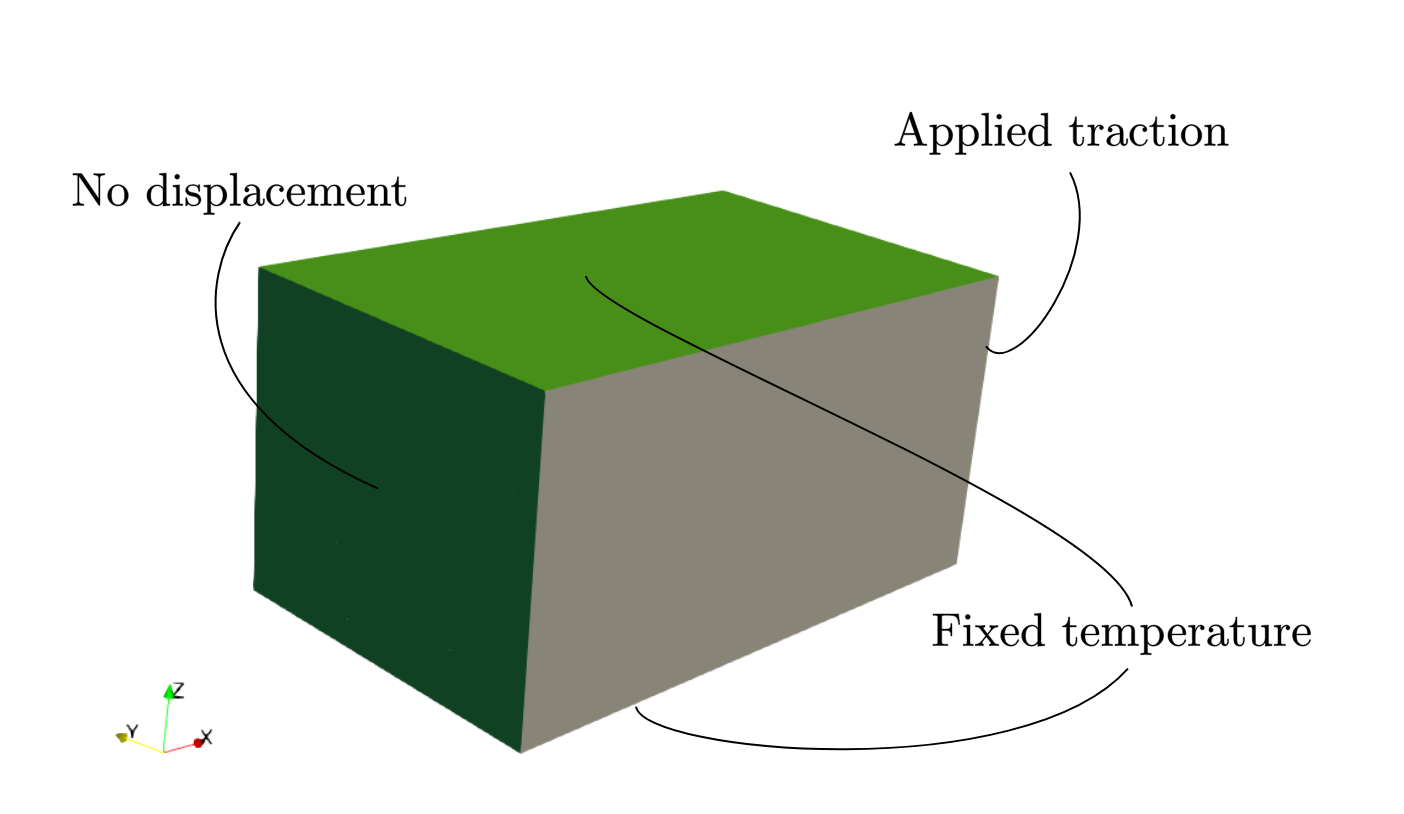}
    \caption{Boundary conditions}
\end{subfigure}
\caption{Geometry and boundary conditions for the numerical examples where only Dirichlet boundary conditions and non-zero Neumann boundary conditions are given}
\label{fig:geomBC}
\end{figure}

A 3x2 configuration with a spacing of 5 mm between the cells is chosen for the battery pack. The cells are 70 mm in height with a diameter of 21 mm. The pack is cooled at the top and bottom assuming a constant fixed temperature. Structurally the pack is analogous to a cantilever beam with one side fixed and the other subjected to a traction load of 15 N/mm$^2$. The boundary conditions and the geometry are presented in Fig.\ref{fig:geomBC}. 

The battery pack is discretized with 48x32x40 linear hexahedra elements leading to a total of 265,188 degrees of freedom. The material properties for the study are given in Table \ref{tab:mat_prop}.

\begin{table}[h]
\begin{center}
\begin{minipage}{\textwidth}
\caption{Properties for the thermo-mechanical analysis}\label{tab:mat_prop}
\begin{tabular*}{\textwidth}{@{\extracolsep{\fill}}lcccccc@{\extracolsep{\fill}}}
\toprule
Property & Aluminum & Cell \\
\midrule
Thermal conductivity $\kappa$ [W/m.K] & 220.0  & 1.25 \\
Young's modulus $E$ [GPa]  & 68.0  & 1.5\\
Poisson's ratio $\nu$ [-] & 0.32  & 0.2 \\
Thermal expansion $\alpha$ [K$^{-1}$] & 21x10$^{-6}$  & 10x10$^{-6}$ \\
\midrule
\end{tabular*}
\end{minipage}
\end{center}
\end{table}

\subsection{Optimization results}
\begin{figure}[!h]
\centering
\begin{subfigure}{0.49\textwidth}
    \includegraphics[width=\textwidth]{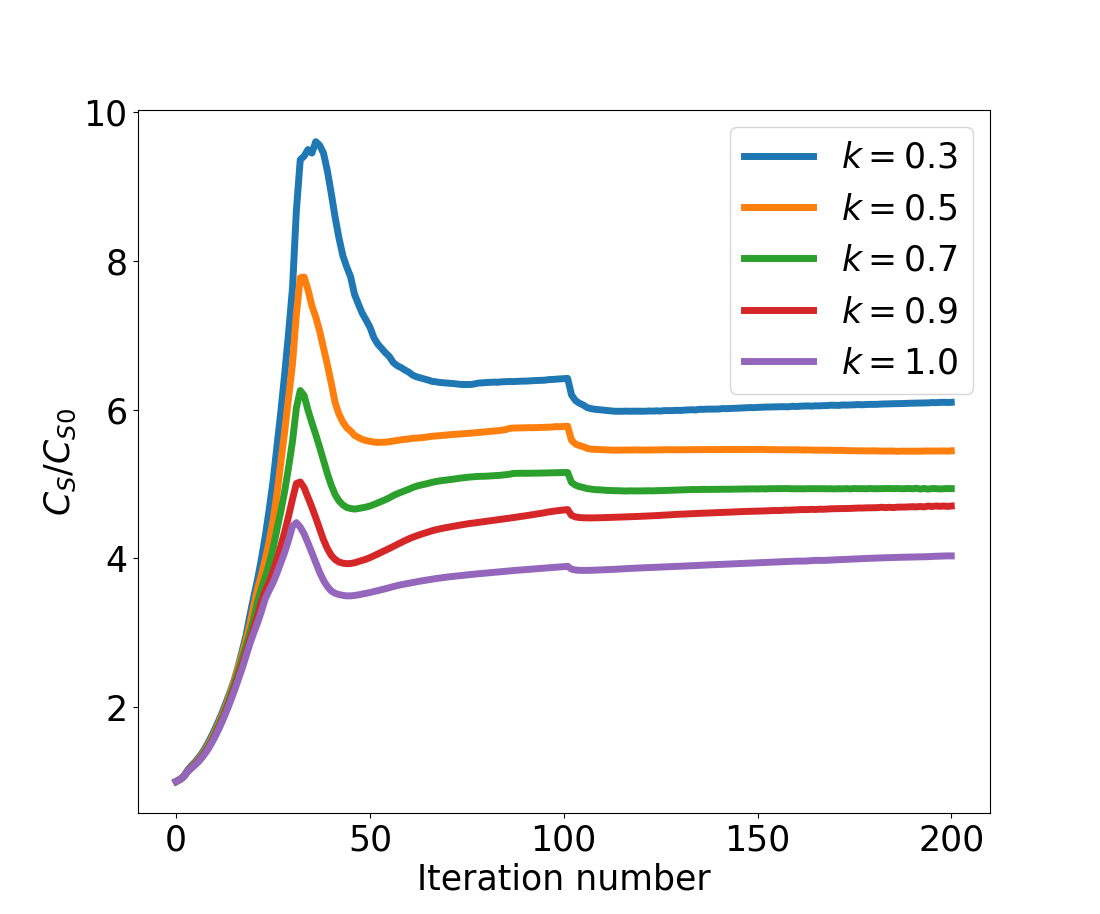}
    \caption{Structural compliance convergence history}
\end{subfigure}
\hfill
\begin{subfigure}{0.49\textwidth}
    \includegraphics[width=\textwidth]{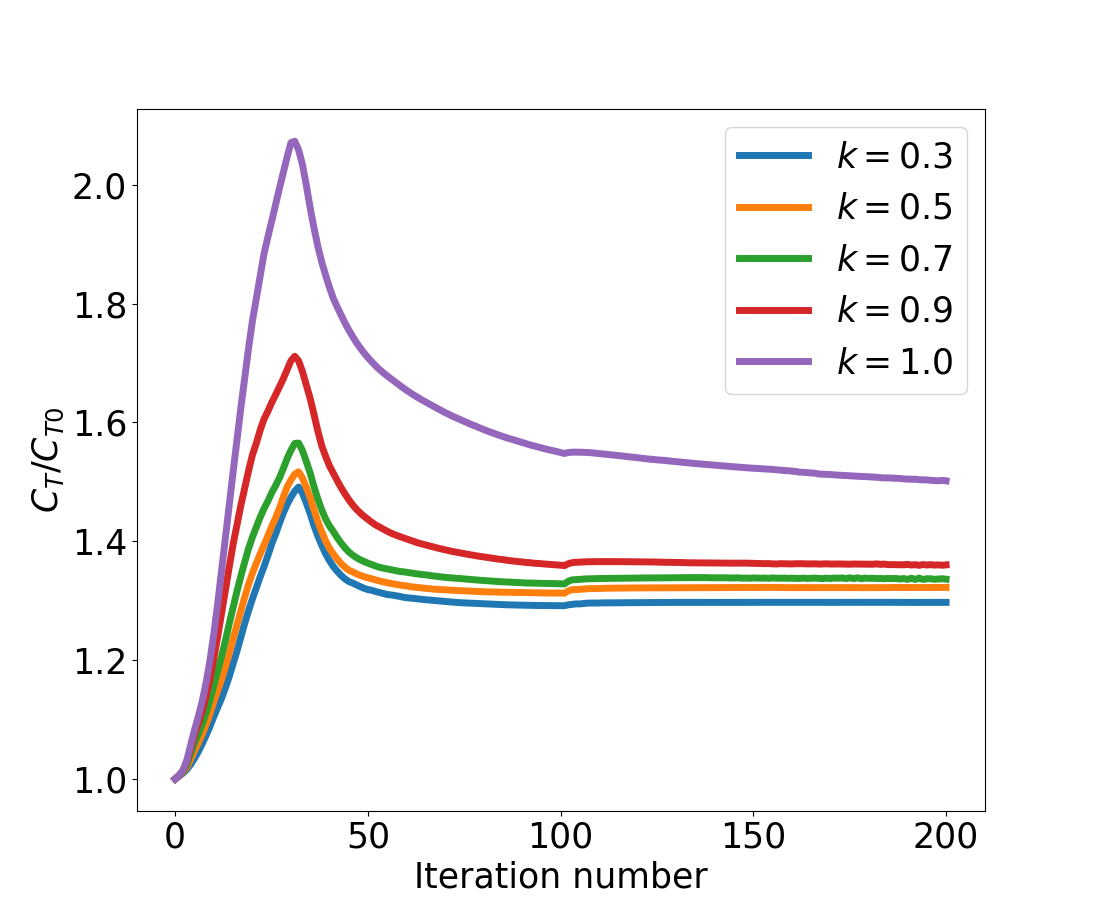}
    \caption{Thermal compliance convergence history}
\end{subfigure}
\caption{Convergence history of structural and thermal compliances for different values of $k$. Note that the jump at iteration 100 is due to a reduction of the move limit value}
\label{fig:Convergence}
\end{figure}

\begin{figure}[!h]
\centering
\begin{subfigure}{0.32\textwidth}
    \includegraphics[width=\textwidth]{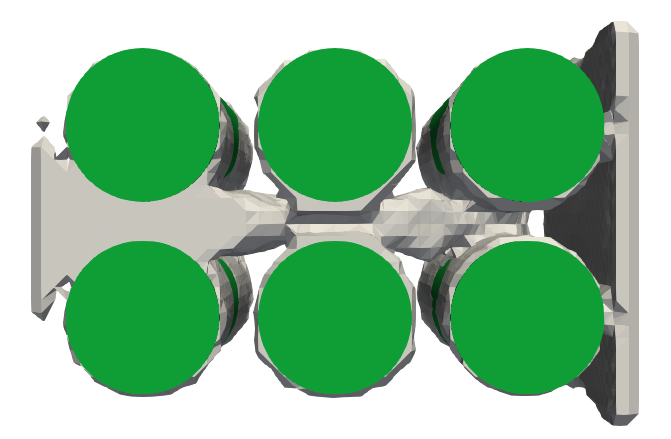}
    \caption{top view $k$ = 0.3}
\end{subfigure}
\hfill
\begin{subfigure}{0.32\textwidth}
    \includegraphics[width=\textwidth]{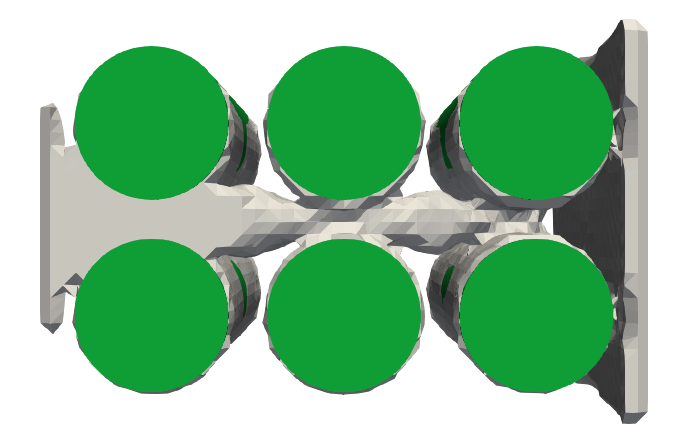}
    \caption{top view $k$ = 0.5}
\end{subfigure}
\hfill
\begin{subfigure}{0.32\textwidth}
    \includegraphics[width=\textwidth]{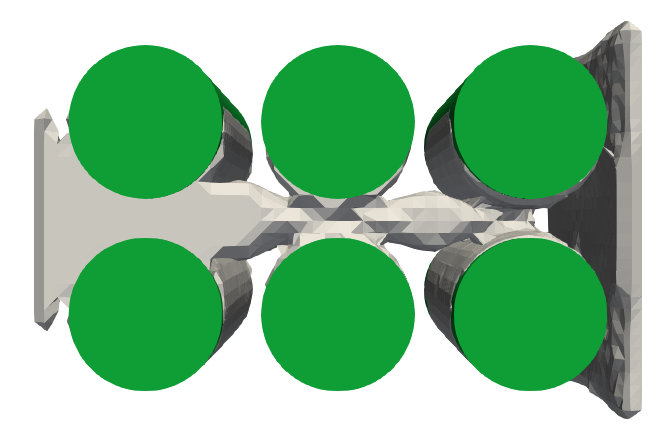}
    \caption{top view $k$ = 1.0}
\end{subfigure}
\begin{subfigure}{0.32\textwidth}
    \includegraphics[width=\textwidth]{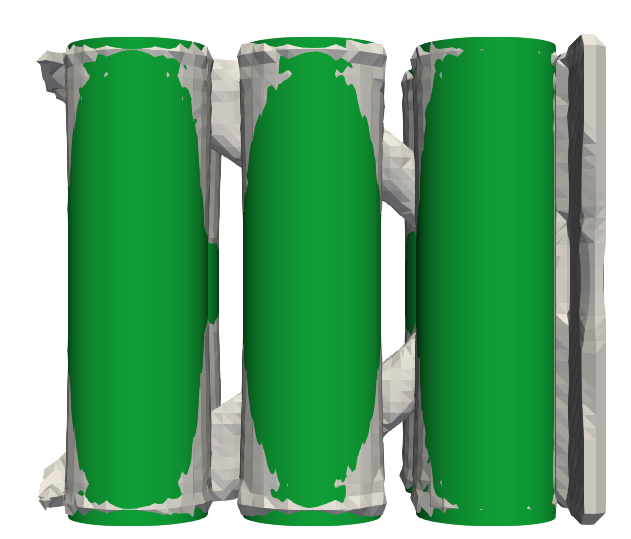}
    \caption{front view $k$ = 0.3}
\end{subfigure}
\hfill
\begin{subfigure}{0.32\textwidth}
    \includegraphics[width=\textwidth]{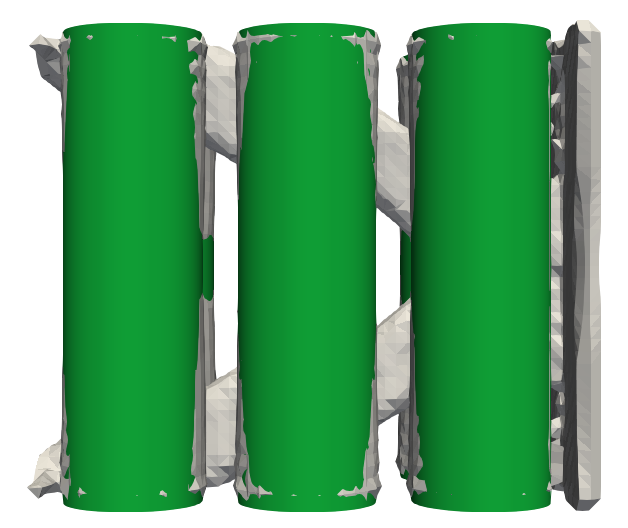}
    \caption{front view $k$ = 0.5}
\end{subfigure}
\hfill
\begin{subfigure}{0.32\textwidth}
    \includegraphics[width=\textwidth]{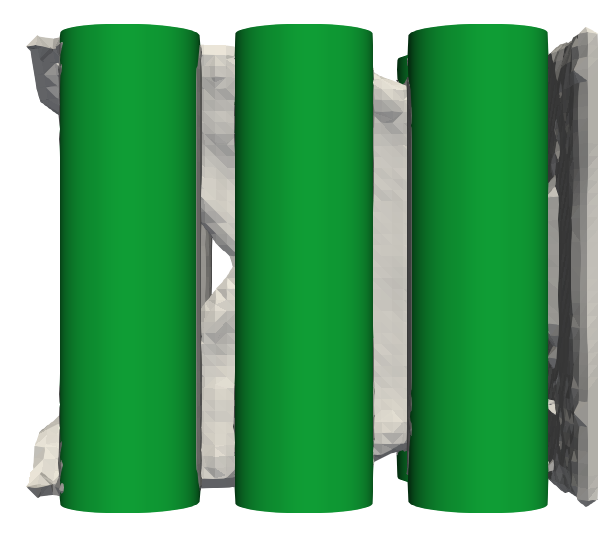}
    \caption{front view $k$ = 1.0}
\end{subfigure}
\begin{subfigure}{0.32\textwidth}
    \centering
    \includegraphics[width=.9\textwidth]{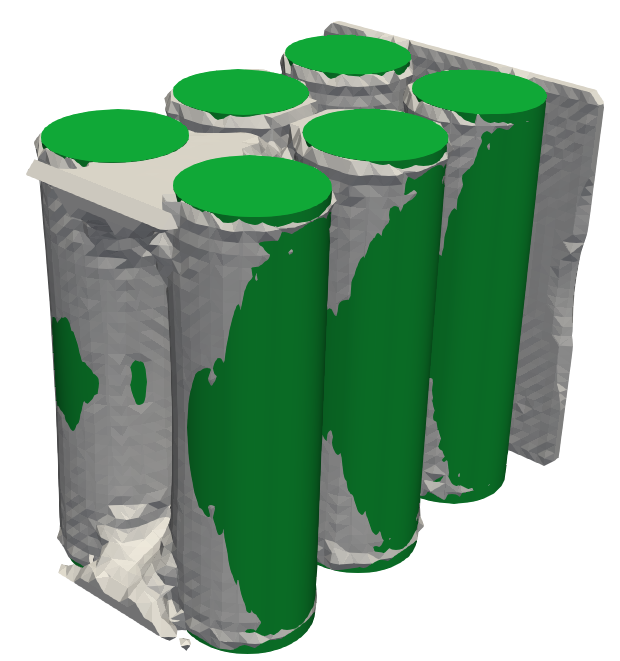}
    \caption{isometric view $k$ = 0.3}
\end{subfigure}
\hfill
\begin{subfigure}{0.32\textwidth}
    \centering
    \includegraphics[width=.9\textwidth]{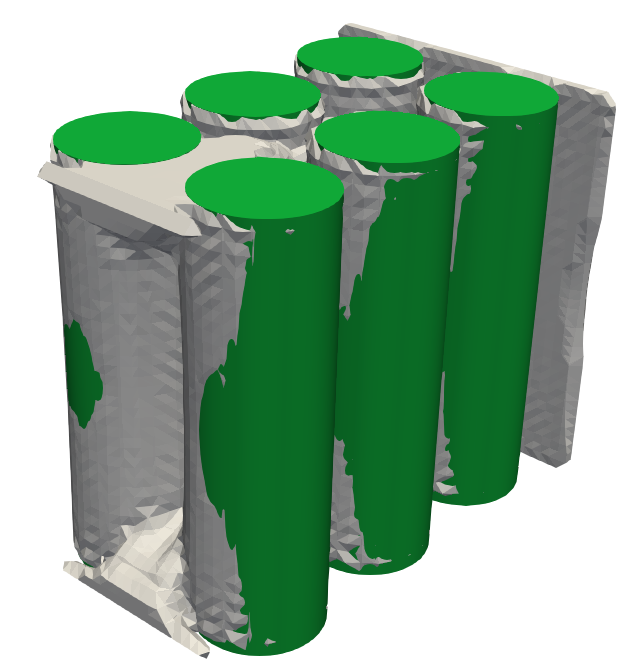}
    \caption{isometric view $k$ = 0.5}
\end{subfigure}
\hfill
\begin{subfigure}{0.32\textwidth}
    \centering
    \includegraphics[width=.9\textwidth]{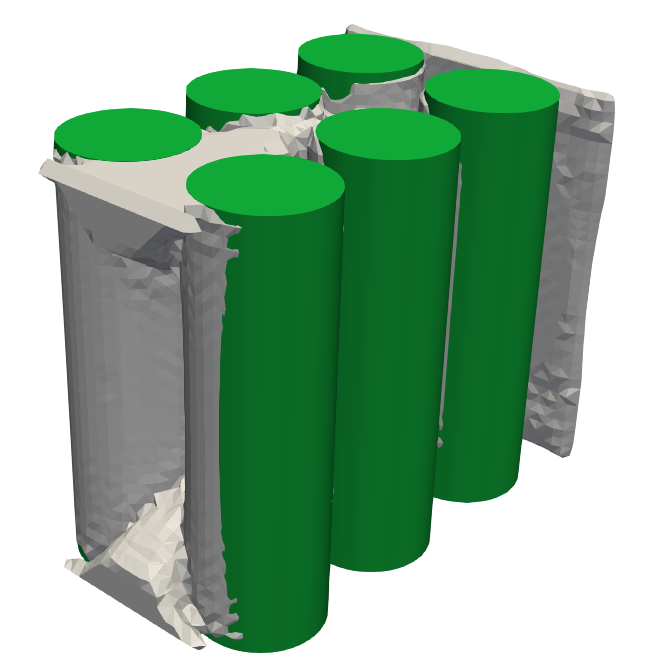}
    \caption{isometric view $k$ = 1.0}
\end{subfigure}
\caption{Optimization results for selected cases}
\label{fig:OptimizationResults}
\end{figure}

\begin{figure}[!h]
\centering
\begin{subfigure}{0.49\textwidth}
    \includegraphics[width=\textwidth]{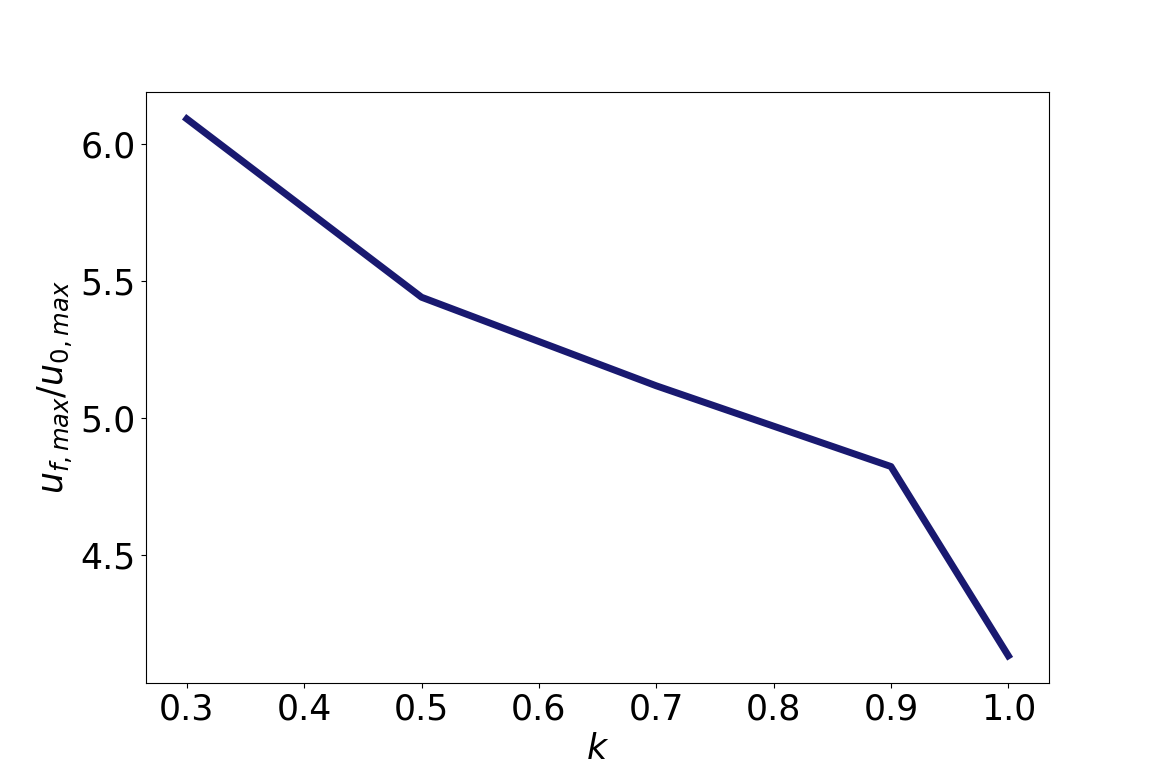}
    \caption{Maximum displacement}
\end{subfigure}
\hfill
\begin{subfigure}{0.49\textwidth}
    \includegraphics[width=\textwidth,  trim={0.8cm 0 0.8cm 0}]{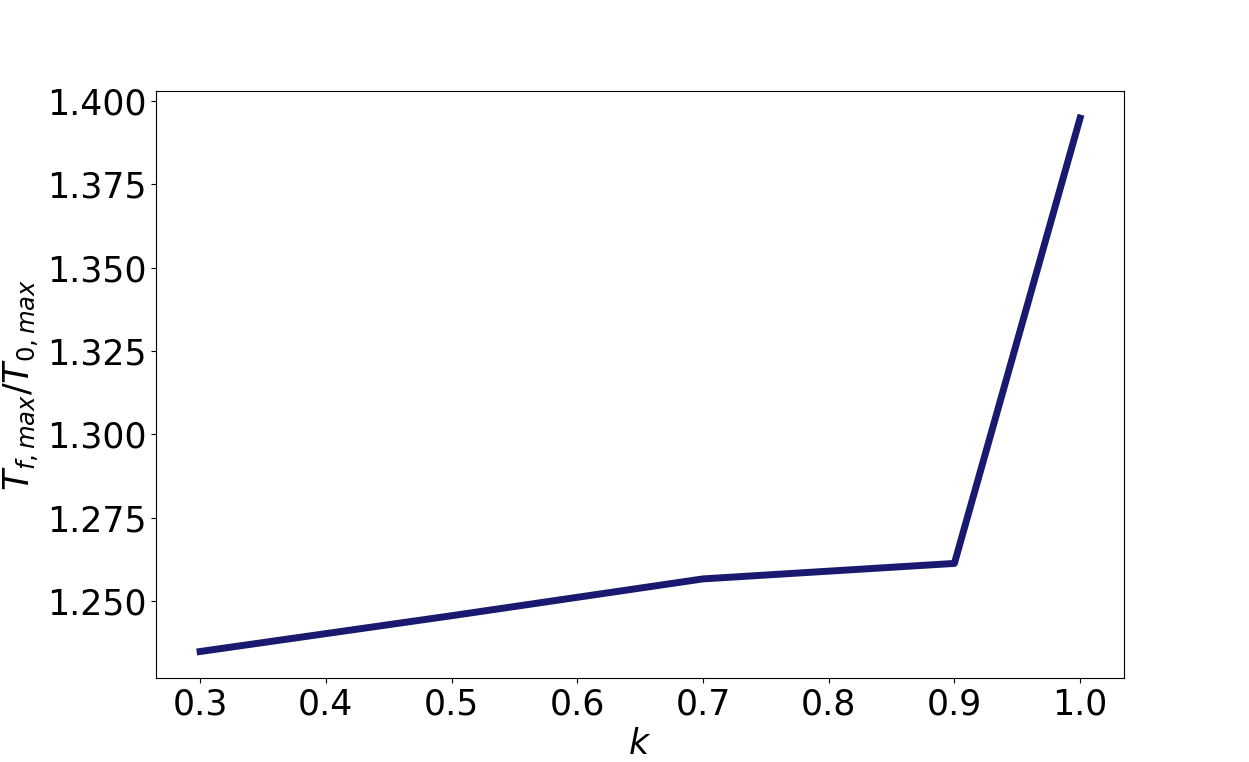}
    \caption{Maximum temperature}
\end{subfigure}
\caption{Maximum displacement and temperature of the final design normalized by the maximum state of the initial design as a function of $k$}
\label{fig:ParetoStates}
\end{figure}

The following weights are investigated $k =$ 1.0, 0.9, 0.7, 0.5, 0.3. The convergence history for the thermal and structural compliance for the 5 cases is given in Fig.\ref{fig:Convergence}. First, the volume constraint is met, and then the objective is gradually minimized until it reaches the minimum bound. The jump at iteration 100 is caused by a modification in the move limit value, which limits the amount of change allowed from one iteration to the next. The optimization results for 3 selected cases are presented in Fig.\ref{fig:OptimizationResults}. As the value of $k$ is decreased, we observe a shift in the topology between the cells - from a stiffener-like structure to one consisting of elements that resemble beams. Additionally, more material is allocated towards the ends of the batteries to facilitate their connection to the heat sinks. Finally, the maximum displacement and maximum temperature of the final design as a function of the weights given to each compliance are given in Fig.\ref{fig:ParetoStates}. Once thermal compliance is taken into account, a significant reduction in the normalized maximum temperature is observed. However, as the value of $k$ continues to decrease, the decrease in maximum temperature becomes less pronounced. This trend is similar to what is observed in the convergence plot shown in Fig.\ref{fig:Convergence}.

\subsection{Transient response of the battery pack}

\begin{figure}[!h]
      \centering
      \includegraphics[width=\textwidth]{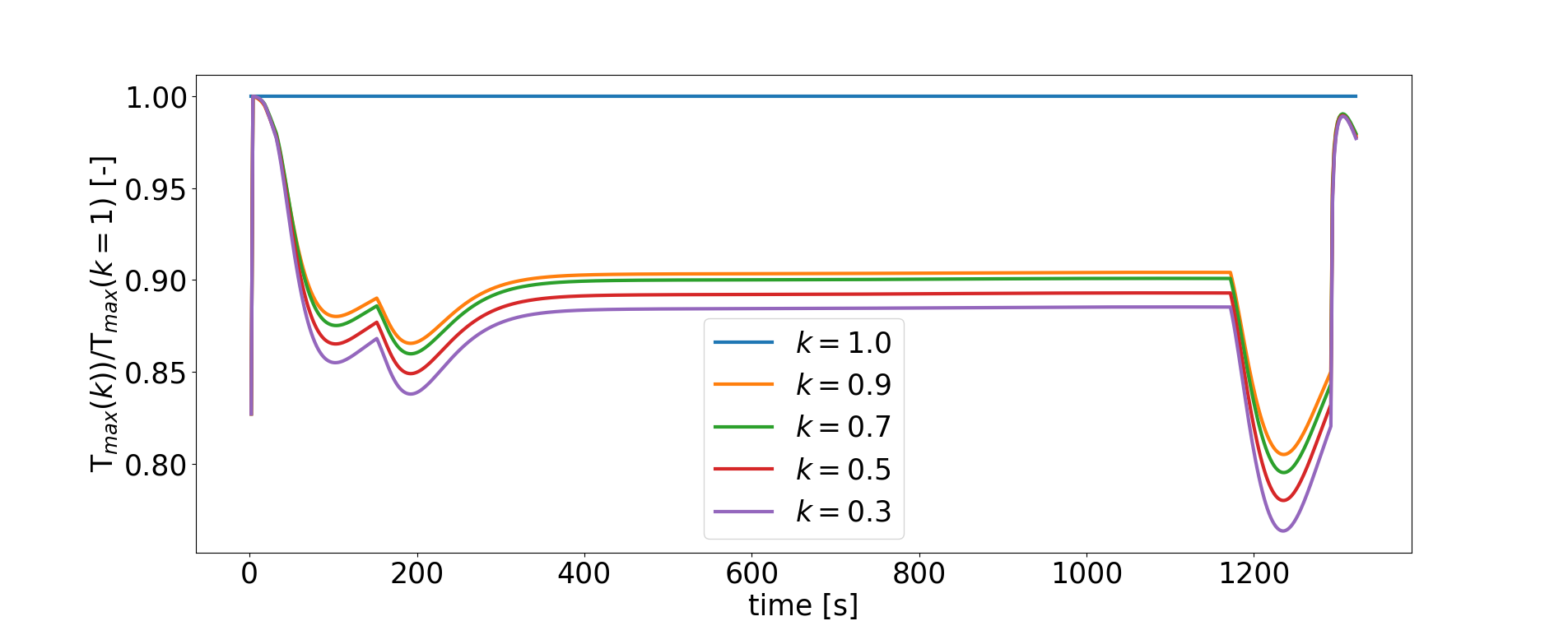}
   \caption{Maximum temperature time history normalized by the maximum temperature time history for $k=1$}
   \label{fig:maxTtime}
\end{figure}

The effective volumetric heat capacity $\rho c_p$ is taken as 1,767,574 J.m$^{-3}$.K$^{-1}$ and 2,430,000 J.m$^{-3}$.K$^{-1}$ for the cell and the pack, respectively \cite{Chen2020, o2022thermal, steinhardt2022meta}. The optimized topologies are analyzed using a transient heat conduction model implemented using FEniCS. An implicit Euler scheme that is unconditionally stable is chosen for the time discretization. The temperature time history and temperature distribution are presented in Fig.\ref{fig:maxTtime}. The maximum temperatures during take-off and landing do not appear to be significantly affected by thermal compliance considerations. These phases are short and the steady state is not reached. The overall maximum temperature does not seem to be much reduced when considering the topology optimized for steady-state considerations. Consequently, it seems worthwhile to consider transient heat conduction topology optimization specifically for takeoff and landing or for the entire flight but that is beyond the scope of this study. On the other hand, in the cruise segment where steady-state is reached, considering thermal compliance in the optimization formulation leads to a reduction of more than 10 \%  of the maximum temperature. Such a reduction could mitigate battery degradation. More research is needed to quantify this effect.    

\section{Conclusions}

In this study, we presented a methodology for conducting thermo-mechanical level-set topology optimization of a load-carrying battery pack in an eVTOL and employed an electrochemical model to predict heat generation from the batteries based on a given flight profile. A multi-objective function that considered both thermal and structural compliance was minimized. We demonstrated the approach for several cases and analyzed them using a transient heat conduction model. The methodology can be adapted and reused for different power profiles, boundary conditions, and pack configurations. 
In the future, we plan to expand our research to encompass various power profiles and boundary conditions, examine degradation and different battery chemistries, and increase the scale of the thermo-mechanical model.

\section*{Acknowledgements}

The authors acknowledge the support from NASA under award No. 80NSSC21M0070.



\bibliography{AeroBest_paper}

\begin{thebibliography}{20}
\providecommand{\natexlab}[1]{#1}
\providecommand{\url}[1]{\texttt{#1}}
\expandafter\ifx\csname urlstyle\endcsname\relax
  \providecommand{\doi}[1]{doi: #1}\else
  \providecommand{\doi}{doi: \begingroup \urlstyle{rm}\Url}\fi

\bibitem[Doyle et~al.(1993)Doyle, Fuller, and Newman]{Doyle1993}
M.~Doyle, T.~F. Fuller, and J.~Newman.
\newblock {Modeling of galvanostatic charge and discharge of the
  lithium/polymer/insertion cell}.
\newblock \emph{Journal of the Electrochemical society}, 140\penalty0
  (6):\penalty0 1526--1533, 1993.
\newblock \doi{10.1149/1.2221597}.

\bibitem[Ai et~al.(2019)Ai, Kraft, Sturm, Jossen, and Wu]{Ai2019}
W.~Ai, L.~Kraft, J.~Sturm, A.~Jossen, and B.~Wu.
\newblock Electrochemical thermal-mechanical modelling of stress inhomogeneity
  in lithium-ion pouch cells.
\newblock \emph{Journal of The Electrochemical Society}, 167\penalty0
  (1):\penalty0 013512, 2019.
\newblock \doi{10.1149/2.0122001JES}.

\bibitem[Timms et~al.(2021)Timms, Marquis, Sulzer, Please, and
  Chapman]{Timms2021}
R.~Timms, S.~G. Marquis, V.~Sulzer, C.~P. Please, and S.~J. Chapman.
\newblock {Asymptotic Reduction of a Lithium-ion Pouch Cell Model}.
\newblock \emph{SIAM Journal on Applied Mathematics}, 81\penalty0 (3):\penalty0
  765--788, 2021.
\newblock \doi{10.1137/20M1336898}.

\bibitem[Deshpande et~al.(2012)Deshpande, Verbrugge, Cheng, Wang, and
  Liu]{Deshpande2012}
R.~Deshpande, M.~Verbrugge, Y.-T. Cheng, J.~Wang, and P.~Liu.
\newblock Battery cycle life prediction with coupled chemical degradation and
  fatigue mechanics.
\newblock \emph{Journal of the Electrochemical Society}, 159\penalty0
  (10):\penalty0 A1730, 2012.
\newblock \doi{10.1149/2.049210jes}.

\bibitem[Sulzer et~al.(2021)Sulzer, Marquis, Timms, Robinson, and
  Chapman]{Sulzer2021}
V.~Sulzer, S.~G. Marquis, R.~Timms, M.~Robinson, and S.~J. Chapman.
\newblock {Python Battery Mathematical Modelling (PyBaMM)}.
\newblock \emph{Journal of Open Research Software}, 9\penalty0 (1):\penalty0
  14, 2021.
\newblock \doi{10.5334/jors.309}.

\bibitem[Andersson et~al.(2019)Andersson, Gillis, Horn, Rawlings, and
  Diehl]{Andersson2019}
J.~A.~E. Andersson, J.~Gillis, G.~Horn, J.~B. Rawlings, and M.~Diehl.
\newblock {CasADi -- A software framework for nonlinear optimization and
  optimal control}.
\newblock \emph{Mathematical Programming Computation}, 11\penalty0
  (1):\penalty0 1--36, 2019.
\newblock \doi{10.1007/s12532-018-0139-4}.

\bibitem[Harris et~al.(2020)Harris, Millman, van~der Walt, Gommers, Virtanen,
  Cournapeau, Wieser, Taylor, Berg, Smith, et~al.]{Harris2020}
C.~R. Harris, K.~J. Millman, S.~J. van~der Walt, R.~Gommers, P.~Virtanen,
  D.~Cournapeau, E.~Wieser, J.~Taylor, S.~Berg, N.~J. Smith, et~al.
\newblock {Array programming with NumPy}.
\newblock \emph{Nature}, 585\penalty0 (7825):\penalty0 357--362, 2020.
\newblock \doi{10.1038/s41586-020-2649-2}.

\bibitem[Aln{\ae}s et~al.(2015)Aln{\ae}s, Blechta, Hake, Johansson, Kehlet,
  Logg, Richardson, Ring, Rognes, and Wells]{alnaes2015fenics}
M.~Aln{\ae}s, J.~Blechta, J.~Hake, A.~Johansson, B.~Kehlet, A.~Logg,
  C.~Richardson, J.~Ring, M.~E. Rognes, and G.~N. Wells.
\newblock The fenics project version 1.5.
\newblock \emph{Archive of Numerical Software}, 3\penalty0 (100), 2015.

\bibitem[Logg and Wells(2010)]{logg2010dolfin}
A.~Logg and G.~N. Wells.
\newblock Dolfin: Automated finite element computing.
\newblock \emph{ACM Transactions on Mathematical Software (TOMS)}, 37\penalty0
  (2):\penalty0 1--28, 2010.

\bibitem[Logg et~al.(2012)Logg, Mardal, and Wells]{logg2012automated}
A.~Logg, K.-A. Mardal, and G.~Wells.
\newblock \emph{Automated solution of differential equations by the finite
  element method: The FEniCS book}, volume~84.
\newblock Springer Science \& Business Media, 2012.

\bibitem[Dunning and Kim(2015)]{dunning2015introducing}
P.~D. Dunning and H.~A. Kim.
\newblock Introducing the sequential linear programming level-set method for
  topology optimization.
\newblock \emph{Structural and Multidisciplinary Optimization}, 51:\penalty0
  631--643, 2015.

\bibitem[Sivapuram et~al.(2016)Sivapuram, Dunning, and
  Kim]{sivapuram2016simultaneous}
R.~Sivapuram, P.~D. Dunning, and H.~A. Kim.
\newblock Simultaneous material and structural optimization by multiscale
  topology optimization.
\newblock \emph{Structural and multidisciplinary optimization}, 54:\penalty0
  1267--1281, 2016.

\bibitem[Kambampati and Jauregui(2023)]{kambampati_sandilya_2023_7613753}
S.~Kambampati and C.~M. Jauregui.
\newblock Paralesto, Feb. 2023.
\newblock URL \url{https://doi.org/10.5281/zenodo.7613753}.

\bibitem[Hyun et~al.(2022)Hyun, Jauregui, Kim, and
  Neofytou]{hyun2022development}
J.~Hyun, C.~Jauregui, H.~A. Kim, and A.~Neofytou.
\newblock On development of an accessible and non-intrusive level-set topology
  optimization framework via the discrete adjoint method.
\newblock In \emph{AIAA SCITECH 2022 Forum}, page 2548, 2022.

\bibitem[Guibert et~al.(2023)Guibert, Hyun, Neofytou, and
  Kim]{guibert2023implementation}
A.~T. Guibert, J.~Hyun, A.~Neofytou, and H.~A. Kim.
\newblock Implementation of a plug-and-play reusable level-set topology
  optimization framework via comsol multiphysics.
\newblock In \emph{AIAA SCITECH 2023 Forum}, page 1675, 2023.

\bibitem[Chen et~al.(2020)Chen, Brosa~Planella, O'Regan, Gastol, Widanage, and
  Kendrick]{Chen2020}
C.-H. Chen, F.~Brosa~Planella, K.~O'Regan, D.~Gastol, W.~D. Widanage, and
  E.~Kendrick.
\newblock {Development of Experimental Techniques for Parameterization of
  Multi-scale Lithium-ion Battery Models}.
\newblock \emph{Journal of The Electrochemical Society}, 167\penalty0
  (8):\penalty0 080534, 2020.
\newblock \doi{10.1149/1945-7111/ab9050}.

\bibitem[Yang et~al.(2021)Yang, Liu, Ge, Rountree, and
  Wang]{yang2021challenges}
X.-G. Yang, T.~Liu, S.~Ge, E.~Rountree, and C.-Y. Wang.
\newblock Challenges and key requirements of batteries for electric vertical
  takeoff and landing aircraft.
\newblock \emph{Joule}, 5\penalty0 (7):\penalty0 1644--1659, 2021.

\bibitem[Batemo(2022)]{batemo_2022}
Batemo.
\newblock Lg chem inr21700-m50l, Oct 2022.
\newblock URL
  \url{https://www.batemo.de/products/batemo-cell-library/lg-chem-inr21700-m50l/}.

\bibitem[O'Regan et~al.(2022)O'Regan, Planella, Widanage, and
  Kendrick]{o2022thermal}
K.~O'Regan, F.~B. Planella, W.~D. Widanage, and E.~Kendrick.
\newblock Thermal-electrochemical parameters of a high energy lithium-ion
  cylindrical battery.
\newblock \emph{Electrochimica Acta}, 425:\penalty0 140700, 2022.

\bibitem[Steinhardt et~al.(2022)Steinhardt, Barreras, Ruan, Wu, Offer, and
  Jossen]{steinhardt2022meta}
M.~Steinhardt, J.~V. Barreras, H.~Ruan, B.~Wu, G.~J. Offer, and A.~Jossen.
\newblock Meta-analysis of experimental results for heat capacity and thermal
  conductivity in lithium-ion batteries: A critical review.
\newblock \emph{Journal of Power Sources}, 522:\penalty0 230829, 2022.

\end{thebibliography}

\end{document}